\documentclass[12pt]{amsart}
\usepackage{amsmath}
\usepackage{amsfonts}
\usepackage{amssymb}
\usepackage{latexsym}
\usepackage[OT2,OT1]{fontenc}
\newcommand\cyr{%
\renewcommand\rmdefault{wncyr}%
\renewcommand\sfdefault{wncyss}%
\renewcommand\encodingdefault{OT2}%
\normalfont
\selectfont}
\DeclareTextFontCommand{\textcyr}{\cyr} 
\def\Om{\Omega}

\def\om{\omega}
\def\a{\alpha}
\def\b{\beta}
\def\g{\gamma}
\def\de{\delta}
\def\var{\varepsilon}

\def\R{{\Bbb R}}

\def\pa{\partial}

\def\Cal{\mathcal}
\def\Lap{\bigtriangleup}
\def\grad{\nabla}
\def\deb{\rightharpoonup}
\def\blackbox{\unskip\kern 6pt\penalty 500%
\raise -1pt\hbox{\vrule\vbox to 8pt{\hrule width 6pt\vfill\hrule}\vrule}}

\newtheorem{thm}{Theorem}

\begin{document}

\title[Particles in a Stokes flow] {The Mean-Field Limit for Solid 
Particles\\ in a Navier-Stokes Flow}

\author{Laurent Desvillettes}
\address[L.D.]{Ecole Normale Sup\'erieure de Cachan\\
CMLA, 61, Av. du Pdt. Wilson\\ F94235 Cachan Cedex}
\email{desville@cmla.ens-cachan.fr}

\author{Fran\c cois Golse}
\address[F.G.]{Laboratoire J.-L. Lions\\
Universit\'e Pierre-et-Marie Curie\\
Bo\^\i te courrier 187\\
75252 Paris Cedex 05\\
\& Centre de Math\'ematiques Laurent Schwartz\\
Ecole Polytechnique, F91128 Palaiseau Cedex}
\email{golse@math.polytechnique.fr}

\author{Valeria Ricci}
\address[V.R.]{ Dipartimento di Metodi e Modelli Matematici\\
Universit\`a di Palermo\\ 
Viale delle Scienze Edificio 8\\
I90128 Palermo}
\email{ricci@unipa.it}

\begin{abstract}
We propose a mathematical derivation of Brinkman's force for a cloud of
particles immersed in an incompressible viscous fluid. Specifically, we
consider the Stokes or steady Navier-Stokes equations in a bounded
domain $\Om\subset\R^3$ for the velocity field $u$ of an incompressible
fluid with kinematic viscosity $\nu$ and density $1$. Brinkman's force
consists of a source term $6\pi\nu j$ where $j$ is the current density
of the particles, and of a friction term $6\pi\nu\rho u$ where $\rho$ is
the number density of particles. These additional terms in the motion
equation for the fluid are obtained from 
the Stokes or steady Navier-Stokes equations set in $\Om$ minus
the disjoint union of $N$ balls of radius $\var=1/N$ in the large
$N$ limit with no-slip boundary condition. The
number density $\rho$ and current density $j$ are obtained from
the limiting phase space empirical measure 
$\tfrac1N\sum_{1\le k\le N}\de_{x_k,v_k}$, where $x_k$ is the
center of the $k$-th ball and $v_k$ its instantaneous velocity.
This can be seen as a generalization of Allaire's result in
[Arch. Rational Mech. Analysis \textbf{113} (1991) 209--259]
who considered the case of periodically distributed $x_k$s
with $v_k=0$, and our proof is based on slightly simpler though
similar homogenization
arguments. Similar equations are used for describing the fluid 
phase in various models for sprays.
\end{abstract}

\maketitle



\section{Introduction} \label{intro}

The subject matter of this paper is the derivation of macroscopic models
for the dynamics of large systems of solid particles or liquid droplets
immersed in a viscous fluid (liquid or gas). Specifically, we are concerned 
with the collective effect of the friction force exerted on each particle as a 
result of the viscosity of the fluid together with a no-slip condition at the 
surface of each particle. This type of fluid/solid interaction is relevant in 
several different physical contexts.

A first example is provided by the sedimentation of solid particles in a
viscous incompressible fluid (say, a liquid), typically under the effect of 
gravity. There is a huge literature on this subject; we shall only mention
a few, such as \cite{Batch}, \cite{Feuille}, \cite{CaflLuke}.

Another example is the case of sprays. 
Sprays are complex flows which 
consist of a dispersed phase immersed in some viscous fluid. 

Sprays can be described (Cf. \cite{orourke, jabin}) by systems of coupled 
macroscopic equations (Eulerian-Eulerian modeling) or by the coupling 
of a macroscopic equation and a kinetic equation (Eulerian-Lagrangian
modeling).

We do not claim that the results in this paper provide a complete derivation
of any of these models for sprays, in particular because we do not analyze
the coupling between the particle and the fluid dynamics. Also, our results
apply to steady regimes only, for reasons that will be discussed below.

The present work is only aimed at providing a rigorous derivation of the
Brinkman force created by a cloud of like spherical particles --- we recall
that this force results from the collective effect of the drag exerted on the
particles by the surrounding fluid. In models for sprays, this Brinkman
force would typically be responsible for the coupling between the motion
of the fluid and that of the dispersed phase.

Our approach of this problem is the homogenization method: we more or
less follow earlier works such as \cite{cioramurat} and \cite{allaire} which
only considered periodic distributions of particles. More precisely, the
reference \cite{cioramurat} established the friction term for the Laplace
equation in a periodically perforated domain with homogeneous Dirichlet 
boundary condition. The case of the Stokes or Navier-Stokes equations
was treated in \cite{allaire} by similar arguments.

The discussion in the present paper differs from \cite{allaire} in two ways.
To begin with, only periodic distributions of particles all of which have the
same velocity (which, by Galilean invariance can be taken as $0$) are
considered in \cite{allaire}. In the present paper, we consider clouds of
particles whose phase space empirical measure converges to some
smooth phase space density. Thus, as long as this (mild) assumption 
is verified, the particles considered here can each have their own
instantaneous velocity\footnote{We mention also the paper \cite{hrlv}, where an analogous problem is considered for the Navier-Stokes equations. This paper (as the references therein)  does not
unfortunately contain any detail about the convergence proof.}. Another difference with \cite{allaire} lies in the
method of proof, which may lead to simplifications here and there.
The reference \cite{allaire} closely followed the argument in 
\cite{cioramurat} by truncating the velocity field in the vicinity of each 
particle, an operation that has the disadvantage of leading to velocity
fields that  fail to satisfy the incompressibility condition. In the present
work, the same goal is achieved by removing to the velocity field
some carefully constructed solenoidal boundary layers so that the 
resulting vector field still satisfies the incompressibility condition.
Hence the pressure can be integrated out, thereby leading to 
somewhat easier computations and avoiding painful estimates. Yet, our analysis borrows a lot from \cite{cioramurat}
and \cite{allaire}, especially in the construction of these boundary
layers.

We found it convenient to describe the cloud of particles through its
empirical measure instead of using (marginals of) its $N$-particle
distribution function, as in \cite{Batch}, \cite{Feuille}, \cite{CaflLuke},
\cite{CaflRubi} --- as a matter of fact, most of these references
assume nearly factorized $N$-particle functions, so that both
viewpoints are essentially equivalent.

\section{Presentation of the model and main results} \label{model}

\subsection{Formal derivation of the model}

Consider a system of $N$ identical rigid spheres in a viscous incompressible
fluid with kinematic viscosity $\nu$ and density $\rho_f$. For simplicity, we
assume that the dynamics of the spheres is given, and we seek the collective
effect on the fluid of the drag force on each sphere. We shall make the two
following scaling assumptions:

a) the speed of the spheres is assumed to be small enough, so that the 
quasi-static approximation holds for the fluid motion, and

b) the collective effect of the drag forces exerted on each sphere is of the
same order of magnitude as the external force field driving the fluid.

First, we outline the quasi-static approximation a). Our starting point is the
set of Navier-Stokes equations
\begin{equation}
\label{InitNS}
\begin{aligned}
\pa_tu+u\cdot\grad_xu+\grad_xp&=\nu\Lap_xu+f\,,\quad\nabla_x\cdot u=0\,,
\\
u(t,\cdot)|_{\pa B(x_k(t),r)}&=\dot{x}_k(t)\,,
\end{aligned}
\end{equation}
where $u\equiv u(t,x)\in\R^3$ and $p\equiv p(t,x)$ are respectively the 
velocity and pressure field in the fluid, while $r$ is the radius of the rigid
balls immersed in the fluid and $x_k(t)$ is the position at time $t$ of the
center of the $k$-th ball $B_{x_k(t),r}$. The density of external force per
unit of mass in the fluid is $f\equiv f(t,x)\in\R^3$.

Notice that, in this model, the effect of solid rotation for each particle is
neglected --- together with the amount of torque particles subject to such
solid rotations would exert on the fluid.

Assume that the motion of the spheres occurs at a time scale that is long
compared to the typical time scale of the external force field $f$. In other
words, we postulate the existence of a small parameter $\tau\ll 1$ such 
that 
$$
x_k(t)=X_k(\tau t)\,.
$$
The quasi-static approximation is obtained as follows: defining the slow 
time variable $T=\tau t$ and $u(t,x)=\tau U(\tau t,x)$, the left-hand side 
of the Navier-Stokes equation is rescaled as
$$
\pa_tu+u\cdot\grad_xu-\nu\Lap_xu=
	\tau^2(\pa_TU+U\cdot\grad_xU)-\tau \nu\Lap_xU\,.
$$
Defining
$$
\tau F(T,x)=f\left(\frac{T}\tau,x\right)\,,\,\,
\tau P(T,x)=p\left(\frac{T}\tau,x\right)\,,\quad\hbox{ and }V_k=\frac{dX_k}{dT}
$$
we recast the Navier-Stokes problem above
\begin{equation}
\label{QS-NS}
\begin{aligned}
\tau^2(\pa_TU+U\cdot\grad_xU)+\tau \grad_xP&=\tau (\nu\Lap_xU+F)\,,
	\quad\nabla_x\cdot U=0\,,
\\
U(T,\cdot)|_{\pa B_{X_k(T),r}}&=V_k(T)\,.
\end{aligned}
\end{equation}
Neglecting all terms of order $O(\tau^2)$ in (\ref{QS-NS}), we arrive at the
quasi-static Stokes problem
\begin{equation}
\label{QS-Stokes}
\begin{aligned}
-\nu\Lap_xU+\grad_xP&=F\,,\quad\nabla_x\cdot U=0\,,
\\
U(T,\cdot)|_{\pa B_{X_k(T),r}}&=V_k(T)\,.
\end{aligned}
\end{equation}
Notice that, in the Stokes problem above, $T$ is only a parameter, so that
$X_k(T)$ and $V_k(T)$ can be regarded as independent. In other words,
in the Stokes problem considered below, it will be legitimate, under the
quasi-static approximation, to consider $X_k$ as a constant and yet to
allow $V_k\not=0$.

This accounts for item a) above in the derivation of our model; let us now
discuss item b), namely the collective effect of the drag force exerted on 
the spheres.

We recall that the drag force exerted on a single sphere of radius $r$ 
immersed in a Stokes fluid with kinematic viscosity $\nu$, density $\rho_f$
is 
$$
6\pi\rho_f\nu rV
$$
where $V$ is the relative velocity of the sphere --- relatively to the speed
of the fluid at infinity: see \cite{LL} \S 20.

Hence the collective force field exerted on the fluid by a system of $N$
identical such spheres with prescribed dynamics is of the order of
$$
6\pi\rho_f\nu Nr\langle V\rangle\,,
$$
where $\langle V\rangle$ is the average relative velocity of the spheres.

In the sequel, we assume that the parameters $\nu$ and $\rho_f$ are 
of order $O(1)$, as well as $\langle V\rangle$, but we are interested in   
situations where $r\ll 1$ (small spheres) and $N\gg 1$ (large number
of spheres). In order for the collective effect of the immersed spheres
to be of the same order as that of the driving external force field, we
postulate (without loss of generality) that
\begin{equation}
\label{mean-field}
Nr\simeq\hbox{Const.}>0\,.
\end{equation}
This scaling assumption leads to the mean field approximation listed
above as b).

\subsection{The quasi-static, mean field limit}

Henceforth we use the sphere radius as the small parameter governing 
all limits of interest here, and denote it by $\var>0$ instead of $r$. Thus
we assume that $N\to\infty$, $\var\to 0$ and 
\begin{equation}
\label{Neps}
N\var=1\,.
\end{equation}
We further assume that the fluid and the particles considered here are 
enclosed in a domain $\Omega$ and denote the volume that is left free
for fluid motion by 
$$
\Omega_\var=\Omega\setminus\bigcup_{k=1}^NB_{x_k,\var}\,.
$$
In this setting, the Stokes problem for the velocity field $u_\var$ and the
pressure field $p_\var$ reads
\begin{eqnarray}
\label{prob}
\left\{
\begin{array}{rcc}
-\Lap u_\var+\grad p_\var&=&g, 
\\
\grad\cdot u_\var&=&0,
\end{array} 
\right.
\qquad\hbox{ on }\Om_{\var} .
\end{eqnarray}
Here, the source term $g$ is the ratio of density of external force per unit 
of mass to the kinematic viscosity. This system is supplemented with a 
no-slip boundary conditions for $u$ on the boundary of $\Omega_{\var}$~:
\begin{eqnarray}
\label{bc}
\left\{
\begin{array}{rll}
u|_{\pa B_{x_k,\var}}& = & v_k, \qquad {\hbox{for}} \qquad k=1,..,N,\\
u|_{\pa\Om} &= &0,
\end{array}
\right.
\end{eqnarray}
where $v_k$ is the instantaneous velocity of the (center of mass of the)
$k$-th sphere.

Denote by
\begin{equation}\label{empi}
 F_N(x,v)=\frac{1}{N} \sum_{k=1}^{N}\delta_{x_k,v_k}(x,v)
\end{equation}
the phase space empirical measure of the system of $N$ spheres and by
\begin{equation}\label{empim}
\rho_{N}(x)=\int_{\R^3}  F_N (x,v) \, dv; 
	\qquad j_{N}(x)=\int_{\R^3}  F_N (x,v) \, v\, dv
\end{equation}
its two first moments.

It will be convenient to consider the natural extension of $u_\var$ to $\Om$
defined by
\begin{eqnarray}\label{ext}
\bar{u}_{\var}(x) = \left\{ 
\begin{array}{lll}
u_{\var}(x)  & {\hbox{ if }} & x \in \Om_{\var}, \\
v_k & {\hbox{ if }} &  x \in B_{x_k,\var}, \qquad k= 1,..,N. 
\end{array} \right.
\end{eqnarray}

As recalled above, the Stokes' computation of the friction exerted on an
immersed sphere by the surrounding viscous incompressible fluid involves
the relative velocity of the sphere to the speed of the fluid at infinity. In order
to extend Stokes' analysis to the mean field situation considered here, we
need to assume that the distance between the immersed particles is large
enough compared to their size. Specifically, we assume that
\begin{equation}
\label{sfere}
\inf_{1\le k\not=l\le N}|x_k-x_l|>2r_\var\hbox{ where }r_\var:=\var^{1/3}\,.
\end{equation}
This assumption on the distance between particles is consistent with the
critical scale for the total number of particles discussed in \cite{JabOtto}.
The assumption (\ref{sfere}) allows considering each particle subject
to a drag force given by Stokes' formula independently of other particles.
Obviously, we do not know whether (\ref{sfere}) is preserved under 
particle motion, and this is why only steady situations are considered
here.

Likewise, we assume for simplicity that the fluid and the particles occupy
a smooth bounded domain $\Om\subset\R^3$, and that there is no direct
interaction between the boundary of $\Om$ and any of the immersed 
particles:
\begin{equation}
\label{sfere2}
\inf_{1\le k\le N}\hbox{dist}(x_k,\pa\Om)>r_\var\,.
\end{equation} 

\begin{thm}\label{TH1}
Let $\Om\subset\R^3$ be a smooth bounded domain, and consider a
system of $N$ balls $B_{x_k,\var}$ for $k=1,\ldots,N$ and $\var=1/N$
included in $\Om$ and satisfying conditions (\ref{sfere})-(\ref{sfere2}).
Assume that the empirical measure $F_N$ has uniformly bounded
kinetic energy
$$
\sup_{N\ge 1}\iint_{\Om\times\R^3}\tfrac12|v|^2F_N(x,v)dxdv<\infty
$$
while the macroscopic density and the current converge weakly in
the sense of measures
$$
\rho_{N} \deb \rho, \qquad j_{N} \deb j\hbox{ as }N\to\infty
$$
with $\rho$ and $j$ continuous on $\bar\Om$.

For each $g \in\! (L^2(\Omega))^3$, let $u_{\var}$ be the unique weak 
solution in $(H^1(\Om_{\var}))^3$ of (\ref{prob}), (\ref{bc}), and define
$\bar{u}_\var$ as in (\ref{ext}). Then, $\bar{u}_{\var}$ converges in 
$(L^2(\Om))^3$ to the solution $U$ of 
\begin{eqnarray}
\label{eqli}
\left\{ \begin{array}{rcl}
-\Lap U + \grad \Pi &=& g+6\,\pi \,( j - \rho\, U), 
\\
\grad\cdot U &=& 0,
\\
U|_{\pa\Om} &=& 0
\end{array} \right.
\end{eqnarray}
\end{thm}

As a matter of fact, the same techniques as in the proof of Theorem 1
allow considering the steady Navier-Stokes, instead of Stokes equations.
The starting point in this case is
\begin{eqnarray}
\label{probNS}
\left\{
\begin{array}{rcc}
u_\var\cdot\grad u_\var-\nu\Lap u_\var+\grad p_\var&=&g, 
\\
\grad\cdot u_\var&=&0,
\end{array} 
\right.
\qquad\hbox{ on }\Om_{\var}.
\end{eqnarray}
In writing the system above, we have retained the kinematic viscosity 
$\nu$ instead of absorbing it in the source term as in the linear, Stokes 
case. Hence, unlike in (\ref{prob}), $g$ is the density of external force 
per unit of mass (instead of its ratio to the kinematic viscosity).

The limiting equations in this case are
\begin{eqnarray}
\label{eqliNS}
\left\{ \begin{array}{rcl}
U\cdot\grad U-\nu\Lap U + \grad \Pi &=& g+6\,\pi\nu \,(j - \rho\, U), 
\\
\grad\cdot U &=& 0,
\\
U|_{\pa\Om} &=& 0
\end{array} \right.
\end{eqnarray}

Let us briefly discuss the uniqueness problem for (\ref{eqliNS}). 
By a standard energy argument, one finds that, if $U_1$ and $U_2$
are weak solutions of (\ref{eqliNS}), they must satisfy
$$
\begin{aligned}
6\pi\nu\int_\Om\rho|U_1-U_2|^2dx&+\nu\|\grad(U_1-U_2)\|_{L^2(\Om)}^2
\\
&\le\|\grad U_1\|_{L^2(\Om)}\|U_1-U_2\|_{L^4(\Om)}^2
\\
&+ \|U_2\|_{L^4(\Om)}
\|\grad(U_1-U_2)\|_{L^2(\Om)}\|U_1-U_2\|_{L^4(\Om)}
\end{aligned}
$$
We first recall (see \cite{Ladyz} p. 9) that
$$
\|U_{\var}\|_{L^4(\Om)}^4
	\leq 4\|U_{\var}\|_{L^2(\Om)}\|\grad U_{\var}\|_{L^2(\Om)}^3\,;
$$
together with the Poincar\'e inequality, this entails  
\begin{equation}
\label{H1L4}
\|U_{\var}\|_{L^4(\Om)}^4\leq 4C_P\|\grad U_{\var}\|_{L^2(\Om)}^4
\end{equation}
where $C_P$ denotes the Poincar\'e constant in the domain $\Om$.
Hence
$$
\begin{aligned}
{}&\nu\|\grad(U_1-U_2)\|_{L^2(\Om)}^2
\\
&\qquad\le 2C_P^{1/2}(\|\grad U_1\|_{L^2(\Om)}+\|\grad U_2\|_{L^2(\Om)})
\|\grad(U_1-U_2)\|^{2}_{L^2(\Om)}\,.
\end{aligned}
$$
Therefore, uniqueness holds for (\ref{eqliNS}) if
$$
\nu\ge 2C_P^{1/2}(\|\grad U_1\|_{L^2(\Om)}+\|\grad U_2\|_{L^2(\Om)})\,.
$$
But the usual energy estimates for either of the weak solutions $U_1$
and $U_2$ shows that
$$
\nu\|\grad U_j\|_{L^2(\Om)}
	\le C_P(\|g\|_{L^2(\Om)}+6\pi\nu\|j\|_{L^2(\Om)})\,.
$$
Finally, uniqueness holds for (\ref{eqliNS}) if
$$
\nu^2\ge 4C_P^{3/2}(\|g\|_{L^2(\Om)}+6\pi\nu\|j\|_{L^2(\Om)})
$$
i.e. for $\nu>\nu_0\equiv\nu_0(\|g\|_{L^2(\Om)},\|j\|_{L^2(\Om)},C_P)$.

\begin{thm}\label{TH2}
Under the same assumptions as in Theorem \ref{TH1} and for each
$\nu>\nu_0(\|g\|_{L^2(\Om)},\|j\|_{L^2(\Om)},C_P)$, consider, for
each $g\in L^2(\Om)$ and each $\var=1/N$, a solution $u_\var$ of
the steady Navier-Stokes equations (\ref{probNS}) with the no-slip 
boundary condition (\ref{bc}). Defining its natural extension to $\Om$
to be $\bar u_\var$ as in (\ref{ext}), one has $\bar u_\var\to u$ in
$L^2(\Om)$ as $\var=1/N\to 0$, where $u$ is the unique weak solution
of (\ref{eqliNS}).

\end{thm}

\section{Method of proof} \label{wee}

In this section, we present the strategy for the proofs of Theorems~1 and 2. 

\subsection{Introducing correctors}

We recall that the weak formulation of the Stokes problem (\ref{prob})-(\ref{bc}) 
is
\begin{equation}
\label{fdstokes}
\int_{\Om_{\var}}\grad u_{\var}\cdot\grad Wdx= \int_{\Om_{\var}}g\cdot Wdx\,,
\end{equation}
while  the weak formulation of the Navier-Stokes problem (\ref{probNS}), (\ref{bc}) 
is
\begin{equation}
\label{fdnavstok}
\nu\int_{\Om_{\var}}\grad u_{\var}\cdot\grad Wdx= 
\int_{\Om_{\var}}u_{\var}\otimes u_{\var}:\grad Wdx
+\int_{\Om_{\var}}g\cdot Wdx\,,
\end{equation}
for each test solenoidal vector field $W\in (H^1_0(\Om_{\var}))^3$, i.e. such 
that $\grad\cdot W=0$.

For each $w\in (\Cal D (\Om))^3$ such that $\nabla\cdot w = 0$, we choose
test vector fields of the form
$$
W_{\var}=w-\Cal{B}_{\var}[w]
$$
where $\Cal{B}_{\var}[w]\in (H^1_0(\Om))^3$ satisfies 
$$
\nabla\cdot\Cal{B}_{\var}[w]=0\hbox{ in }\Om\hbox{ and }
\Cal{B}_{\var}[w]\mid_{\bar{B}_{x_k,\var}}= w \mid_{\bar{B}_{x_k,\var}}\,.
$$

Similarly, we approximate the solution by
$$
U_{\var}=\bar{u}_{\var}-\Cal{A}_{\var}
$$
where $\Cal{A}_{\var}\in (H^1_0(\Om))^3$ satisfies
$$
\nabla\cdot\Cal{A}_{\var}=0\hbox{ in }\Om\hbox{ and }
\Cal{A}_{\var}\mid_{\bar{B}_{x_k,\var}}= v_k\,.
$$

Explicit formulas for the fields  $\Cal{A}_{\var}$ and  $\Cal{B}_{\var}$ will
be given at the end of the present section. Notice that, by construction,  
$$
U_{\var} \mid_{B_{x_k,\var}}=W_{\var} \mid_{B_{x_k,\var}}=0\,,\quad
\hbox{ for all }k=1,..,N\,.
$$ 

In addition, the correctors $\Cal{A}_{\var}$ and $\Cal{B}_{\var}$ are chosen
so that
\begin{eqnarray}
\label{B}
\Cal{B}_{\var}[w]\deb 0&\mathrm{in} & (H_{0}^1(\Om))^3\,,
\\
\label{A}
\Cal{A}_{\var}\deb 0 &\mathrm{in} & (H_{0}^1(\Om))^3\,.
\end{eqnarray}

Condition (\ref{B}) implies that 
$$
W_{\var}\deb w\hbox{ in }(H_{0}^1(\Om))^3\hbox{ and }
\nabla W_{\var}\deb\nabla w\hbox{ in }(L^2(\Om))^9\,.
$$
Moreover, (\ref{B}) and (\ref{A}) imply that
\begin{eqnarray*}
\label{B0}
\Cal{B}_{\var}[w]\to 0&\mathrm{in}&(L^p(\Om))^3\,,
\\
\label{A0}
\Cal{A}_{\var}\to 0&\mathrm{in}&(L^p(\Om))^3\,,
\end{eqnarray*}
for each $p\in[1,6)$, by the Rellich-Kondrachov compact embedding theorem,
so that
$$
W_{\var}\to w\hbox{ in }(L^p(\Om))^3\hbox{ for each }p\in[1,6)\,.
$$

Condition (\ref{A}) implies that $U_\var$ and $\bar{u}_\var$ behave similarly
as $\var\to 0$. In the next subsection, we study the asymptotic behavior of 
$U_\var$, which is somewhat simpler to analyze. As we shall see, condition 
(\ref{A}) implies that 
$$
U_{\var}\deb U\hbox{ in }(H_{0}^1(\Om))^3\hbox{ as }\var\to 0
$$ 
for both problems (\ref{prob}) and (\ref{probNS}) with the boundary condition 
(\ref{bc}). Hence
$$
U_{\var} \to U\hbox{ in }(L^p(\Om))^3\hbox{ for }1\le p<6
$$
as $\var\to 0$.

\subsection{Weak convergence of $U_{\var}$} \label{cf}

Here we show that (some subsequence of) $U_{\var}$ converges weakly 
in $(H^1(\Om))^3$ (assuming (\ref{B}) and (\ref{A})), for both problems 
(\ref{prob}) and (\ref{probNS}) with boundary condition (\ref{bc}).

Indeed, for each $k=1,\ldots,N$, one has $U_{\var} \mid_{B_{x_k,\var}}=0$,
so that the weak formulation of the Stokes problem becomes~:
\begin{eqnarray*}
\|\grad U_{\var}\|_{L^2(\Om)}^2
	&=&\int_{\Om}\grad\bar{u}_{\var}:\grad U_{\var}dx
		-\int_{\Om}\grad\Cal{A}_{\var}:\grad U_{\var}dx
\\
&=&\int_{\Om}g\cdot U_{\var}dx
	-\int_{\Om}\grad \Cal{A}_{\var}:\grad U_{\var}dx
\\
&\leq&\|g\|_{L^2(\Om)}\|U_{\var}\|_{L^2(\Om)}
	+\|\grad \Cal{A}_{\var} \|_{L^2(\Om)}\|\grad U_{\var}\|_{L^2(\Om)}  
\end{eqnarray*}
By the Poincar\'e inequality and (\ref{A}), which entails a uniform bound of the
form $\|\grad \Cal{A}_{\var} \|_{L^2(\Om)}<C$, 
we conclude that 
$\|\grad U_{\var}\|_{L^2(\Om)}$ is bounded. Hence there is a subsequence 
such that $U_{\var}\deb U$ in $(H^1(\Om))^3$.

For the Navier-Stokes problem we have similarly:
$$
\begin{aligned}
{}&\nu\|\grad U_{\var}\|_{L^2(\Om)}^2 
=\nu\int_{\Om}\grad \bar{u}_{\var}:\grad U_{\var}dx
-\nu\int_{\Om}\grad \Cal{A}_{\var}:\grad U_{\var}dx
\\
&=\int_{\Om}g\cdot U_{\var}dx
-\int_{\Om_{\var}}(\bar{u}_{\var}\cdot\grad\bar{u}_{\var})\cdot U_{\var}dx
-\nu\int_{\Om}\grad \Cal{A}_{\var}:\grad U_{\var}dx
\\
&=\int_{\Om}g\cdot U_{\var}dx
+\int_{\Om_{\var}}U_{\var}\otimes U_{\var}:\grad U_{\var}dx
-\nu\int_{\Om}\grad \Cal{A}_{\var}:\grad U_{\var}dx
\\
&+\int_{\Om_{\var}}(\Cal{A}_{\var}\otimes U_{\var} 
+U_{\var}\otimes\Cal{A}_{\var}+\Cal{A}_{\var}\otimes\Cal{A}_{\var})
:\grad U_{\var}dx\,.
\end{aligned}
$$
Observe that
$$
\begin{aligned}
\int_{\Om_{\var}}U_{\var}\otimes U_{\var}:\grad U_{\var}dx
&=
\int_{\Om_{\var}}U_{\var}\cdot\left((U_{\var}\cdot\grad)U_{\var}\right)dx
\\
&=
\tfrac12\int_{\Om_{\var}}\grad\left(U_{\var}|U_{\var}|^2\right)dx=0
\end{aligned}
$$
by Green's formula, since $U_\var\big|_{\pa\Om_\var}=0$. Hence 
$$
\begin{aligned}
\nu\|\grad &U_{\var}\|_{L^2(\Om)}^2 \leq\|g\|_{L^2(\Om)}\|U_{\var}\|_{L^2(\Om)}+\nu\|\grad\Cal{A}_{\var}\|_{L^2(\Om)}
\|\grad U_{\var}\|_{L^2(\Om)}    
\\
&+
(\|\Cal{A}_{\var}\otimes U_{\var}\|_{L^2(\Om)}
+\|U_{\var}\otimes\Cal{A}_{\var}\|_{L^2(\Om)}
+\|\Cal{A}_{\var}^{\otimes 2}\|_{L^2(\Om)})\|\grad U_{\var}\|_{L^2(\Om)}
\\
&\leq\|g\|_{L^2(\Om)}\|U_{\var}\|_{L^2(\Om)}
\\
&+(\nu\|\grad\Cal{A}_{\var}\|_{L^2(\Om)}
+2\|\Cal{A}_{\var}\|_{L^4(\Om)}\|U_{\var}\|_{L^4(\Om)}
+\|\Cal{A}_{\var}\|_{L^4(\Om)}^2)\|\grad U_{\var}\|_{L^2(\Om)}\,.
\end{aligned}
$$
Applying inequality (\ref{H1L4}) shows that
\begin{equation}
\begin{aligned}
(\nu-2\sqrt{2}C_P^{1/4}\|\Cal{A}_{\var}\|_{L^4(\Om)})&\|\grad U_{\var}\|_{L^2(\Om)}
\\
&\leq(C_P \|g\|_{L^2(\Om)}+\nu\|\grad \Cal{A}_{\var}\|_{L^2(\Om)}
+\|\Cal{A}_{\var}\|_{L^4(\Om)}^2).
\end{aligned}
\end{equation}
Recall that $\|\Cal{A}_{\var}\|_{L^4(\Om)}\to 0$ as $\var\to 0$ by (\ref{A0}), while
$\|\grad \Cal{A}_{\var}\|_{L^2(\Om)}\le C$ by (\ref{A}). Hence the estimate above 
entails the bound
$$
\|\grad U_\var\|_{L^2(\Om)}\le C\,.
$$

\subsection{Weak formulations on the whole domain}
Next we recast the weak formulations (\ref{fdstokes}) and (\ref{fdnavstok}) in 
terms of $U_\var$: as we shall see, this is somewhat more convenient, at least
in taking the mean field limit.

We first discuss the Stokes problem (\ref{fdstokes}).  Observe that
$$
\int_{\Om}\grad \bar{u}_{\var}:\grad W_{\var}dx= 
\int_{\Om_{\var}} \grad u_{\var}:\grad W_{\var}dx=
\int_{\Om_{\var}} g \cdot W_{\var}dx\,.
$$
Expressing $\bar{u}_\var$ in terms of $U_\var$, one arrives at 
$$
\int_{\Om}\grad U_{\var}:\grad W_{\var}dx 
+\int_{\Om}\grad \Cal{A}_{\var}:\grad W_{\var}dx
=\int_{\Om}g\cdot W_{\var}dx
$$
which, after replacing $W_\var$ with its expression in terms of $w$, leads to
$$
\int_{\Om}\grad U_{\var}:\grad wdx
-\int_{\Om}\grad U_{\var}:\grad \Cal{B}_{\var}dx
+\int_{\Om}\grad \Cal{A}_{\var}:\grad W_{\var}dx
=\int_{\Om}g\cdot W_{\var}dx\,.
$$
Transforming the second integral on the right hand side by Green's formula,
one eventually obtains
\begin{equation}
\label{stokesdeb}
\int_{\Om}\grad U_{\var}:\grad wdx+\int_{\Om}U_{\var}\cdot\Lap\Cal{B}_{\var}dx -\int_{\Om}\Lap \Cal{A}_{\var}\cdot W_{\var}dx=\int_{\Om}g\cdot W_{\var}dx\,.
\end{equation}
Under assumptions (\ref{B})-(\ref{A}) 
\begin{eqnarray}
\int_{\Om}\grad U_{\var}:\grad wdx&\to&\int_{\Om}\grad U:\grad wdx
\\
\int_{\Om}g\cdot W_{\var}dx&\to&\int_{\Om}g\cdot wdx
\end{eqnarray}
as $\var \to 0$. Thus we are left with computing the limit of  
$$
-\int_{\Om}\grad U_{\var}:\grad \Cal{B}_{\var}dx
+
\int_{\Om}\grad\Cal{A}_{\var}:\grad W_{\var}dx
$$ 
or, equivalently, of
$$
\int_{\Om}U_{\var}\cdot\Lap\Cal{B}_{\var}dx
	-\int_{\Om}\Lap\Cal{A}_{\var}\cdot W_{\var}dx\,.
$$

For the Navier-Stokes problem (\ref{fdnavstok}), we follow the same arguments.
First
\begin{equation}
\label{ufo}
\begin{aligned}
{}&\nu\int_{\Om_{\var}}\grad u_{\var}\cdot\grad W_{\var}dx- 
\int_{\Om_{\var}}u_{\var}\otimes u_{\var}:\grad W_{\var}dx
\\
&=\nu\int_{\Om}\grad \bar{u}_{\var}\cdot\grad W_{\var}dx
-\int_{\Om}\bar{u}_{\var}\otimes\bar{u}_{\var}:\grad W_{\var}dx
=\int_{\Om}g\cdot W_{\var}dx
\end{aligned}
\end{equation}
since $v_k$, $k=1,\ldots,N$, are constants and $W_{\var} \mid_{B_{x_k,\var}}=0$.
Making the substitution $\bar{u}_{\var}=U_{\var}+\Cal{A}_{\var}$ in (\ref{ufo}),
one gets
\begin{eqnarray*}
\nu\int_{\Om}\grad U_{\var}:\grad W_{\var}dx
+\nu\int_{\Om}\grad \Cal{A}_{\var}:\grad W_{\var}dx
=
\int_{\Om}U_{\var}\otimes U_{\var}:\grad W_{\var}dx&&
\\+
\int_{\Om}g\cdot W_{\var}dx
+\int_{\Om}(\Cal{A}_{\var}\otimes U_{\var}+U_{\var}\otimes\Cal{A}_{\var}
+\Cal{A}_{\var}\otimes\Cal{A}_{\var}):\grad W_{\var}dx&&
\end{eqnarray*}
and, inserting $W_{\var}=w-\Cal{B}_{\var}$ in the equality above, one eventually
arrives at
\begin{eqnarray}
\label{NS}&&\\
\nonumber
\nu\int_{\Om}\grad U_{\var}:\grad wdx
-\nu\int_{\Om}\grad U_{\var}:\grad \Cal{B}_{\var}dx
+\nu\int_{\Om}\grad \Cal{A}_{\var}:\grad W_{\var}dx&&
\\
\nonumber
=\int_{\Om}U_{\var}\otimes U_{\var}:\grad W_{\var}dx+
\int_{\Om}g\cdot W_{\var}dx&&
\\
+\int_{\Om}(\Cal{A}_{\var}\otimes U_{\var}+U_{\var}\otimes\Cal{A}_{\var}
+\Cal{A}_{\var}\otimes\Cal{A}_{\var}):\grad W_{\var}dx\,.&& 
\nonumber
\end{eqnarray}

Next we pass to the limit as $\var\to 0$; assumptions (\ref{B})-(\ref{A}) imply
that $\Cal{A}_{\var}$ and $U_{\var}$ converge strongly in $L^4(\Om)$, so that
\begin{eqnarray*}
\int_{\Om}(\Cal{A}_{\var}\otimes U_{\var}+U_{\var}\otimes\Cal{A}_{\var}
+\Cal{A}_{\var}\otimes \Cal{A}_{\var}):\grad W_{\var}dx&\to&0\,,
\\
\int_{\Om}U_{\var}\otimes U_{\var}:\grad W_{\var}dx&\to&
\int_{\Om}U\otimes U:\grad wdx\,.
\end{eqnarray*}
Moreover
\begin{eqnarray*}
\int_{\Om}g\cdot W_{\var}dx&\to&\int_{\Om}g\cdot wdx\,,
\\
\nu\int_{\Om}\grad U_{\var}:\grad wdx&\to& 
\nu\int_{\Om}\grad U:\grad wdx\,,
\end{eqnarray*}
so that we are left with the task of computing the limit as $\var\to 0$ of 
$$
-\nu\int_{\Om}\grad U_{\var}:\grad\Cal{B}_{\var}dx
+\nu \int_{\Om}\grad\Cal{A}_{\var}:\grad W_{\var}dx\,.
$$
At this point, we need to specify how the correctors $\Cal{A}_{\var}$ and 
$\Cal{B}_{\var}$ are constructed.

\subsection{Defining the correctors}

Given any smooth function $w$ on $B_{0,s}$ and $r>s$, we designate by
$\Psi_{s,r}[w]$ the solution of the following Stokes problem~:
\begin{eqnarray}
\label{probcor3}
\left\{
\begin{array}{ll}
\Lap\Psi_{s,r}[w]=\grad{\Pi}_{s,r}[w]\,, &x\in B_{0,r}\setminus B_{0,s}\,,
\\
\grad\cdot\Psi_{s,r}[w]=0\,,&
\\
\Psi_{s,r}[w]\mid_{B_{0,s}}=w\,,&
\\
\Psi_{s,r}[w]\mid_{B^{c}_{0,r}}=0\,.&
\\
\end{array}
\right.
\end{eqnarray}
When $s=\var$, $r=r_{\var}=\var^{1/3}$, we define 
$$
\psi_{\var}[w]=\Psi_{\var,\var^{1/3}}[w]\,.
$$
We denote $\pi_{\var}[w]={\Pi}_{\var,\var^{1/3}}[w]$), the pressure field 
associated to $\psi_{\var}[w]$.

With the function $\psi_{\var}[w]$, we define the corrector $\Cal{B}_{\var}$ 
as follows:
\begin{equation}
\label{Bexp}
\Cal{B}_{\var}=\sum_{k=1}^{N}\psi_{\var}[w(\cdot+x_k)](x-x_k)\,.
\end{equation}
\medskip

Whenever $w$ is a constant, i.e. $w(x)=v$, we use the notation
$$
\Phi_{s,r}[v]=\Psi_{s,r}[w]\hbox{ and likewise }\phi_{\var}[v]=\psi_{\var}[w]\,.
$$

With the function $\phi_{\var}[v]$, we define the corrector $\Cal{A}_{\var}$
in the following manner:
\begin{equation}
\label{Axp}
\Cal{A}_{\var}=\sum_{k=1}^{N}\phi_{\var}[v_k](x-x_k)\,.
\end{equation}

The vector fields $\Cal{A}_{\var}$ and $\Cal{B}_{\var}$ so defined are 
obviously solenoidal elements of $(H^1_0(\Om))^3$ that verify the 
conditions 
$$
\Cal{A}_{\var}\mid_{\bar{B}_{x_k,\var}}= v_k\hbox{ and }
\Cal{B}_{\var}[w]\mid_{\bar{B}_{x_k,\var}}=w\mid_{\bar{B}_{x_k,\var}}\,.
$$
In section \ref{esti}, we shall prove that  $\Cal{A}_{\var}$ and $\Cal{B}_{\var}$ 
verify assumptions (\ref{A})-(\ref{B}).

\section{Explicit formulas for the correctors}\label{soluz}

The Stokes equations in an annulus can be solved explicitly; in this section, 
we use these explicit formula to express the correctors $\Cal{A}_{\var}$ and
$\Cal{B}_{\var}$, and to estimate the quantity
$$
\int_{\Om}\grad U_{\var}:\grad \Cal{B}_{\var}dx
	-\int_{\Om}\grad\Cal{A}_{\var}:\grad W_{\var}dx\,.
$$
Occasionally, we will refer to the appendix (section \ref{append}) where a few
standard computations are summarized.

We start with a formula for $\Phi_{1,R}[v]$ --- and hence for $\phi_{\var}[v]$.
For each $x\in \R^3$, denote $r= |x|$, $\om = \frac{x}{|x|}$. Moreover, we
denote $P_{\om}a = (\om\cdot a)\,\om$ is the orthogonal projection on the 
line $\R\,\om$. 

Whenever $1\leq r \leq R$,
\begin{eqnarray} \label{pipo}
\Phi_{1,R}[v] (x)
&=&-\,\bigg[4\,\a(R)\,r^2+2\,\b(R)+\frac{\g(R)}{r}-\frac{\de(R)}{r^3}\bigg]\,
(I-P_{\om})\,v
\nonumber
\\
&&-2\,\bigg[\a(R)\,r^2+\b(R)+\frac{\g(R)}{r}+\frac{\de(R)}{r^3}\bigg]\, 
P_{\om}\,v\,,
\end{eqnarray}
while
\begin{eqnarray*}
\Phi_{1,R}[v] (x)&=&v\hbox{ for }x\in B_{0,1}\,,
\\
\Phi_{1,R}[v] (x)&=&0\hbox{ for }x\in B^c_{0,R}\,.
\end{eqnarray*}

In the formulas above 
\begin{equation}
\label{pippo}
\a(R)=-\frac{3}{8R^3}+O(1/R^4)\,,
\qquad 
\b(R)=\frac{9}{8R}+O(1/R^2)\,,
\end{equation}
while
\begin{equation}
\label{pipppo}
\g(R)=-\tfrac34+O(1/R)\,,\qquad\de(R)=\tfrac14+O(1/R)
\end{equation}
as $R\to+\infty$.

If one replaces the boundary condition at $r=R$ with the condition at infinity
$$
\lim_{|x|\to\infty}\Phi= 0
$$
the solution is 
\begin{equation}
\label{infinito}
\Phi_{1,\infty}[v](x)=\tfrac{1}{4}\,\left(\frac{3}{r}+\frac{1}{r^3}\right)\,(I-P_{\om})\,v
+\tfrac{1}{2}\,\left(\frac{3}{r}-\frac{1}{r^3}\right)\,P_{\om}\, v.
\end{equation}
We denote by $\Pi_{1,\infty}[v]$ the associated pressure.

The following relations hold between the pressure fields $\Pi_{1,R}[v]$ and 
$\Pi_{1,\infty}[v]$~:
\begin{equation}\label{oubli}
\begin{aligned}
\om\cdot\nabla\Phi_{1,R}[v](x)-\Pi_{1,R}[v](x)\,\om
&=
\om\cdot\nabla\Phi_{1,\infty}[v](x)-\Pi_{1,\infty}(x)\,\om 
\\
&-8\,\a(R)\,r \,(I-3P_{\om})\,v+\frac{1}{r^2}\,O\left(\frac{1}{R}\right)
\end{aligned}
\end{equation}
in the limit as $R\to+\infty$.
Finally,
\begin{equation}
\label{forza}
(\om\cdot\nabla\Phi_{1,\infty}[v](x)-\Pi_{1,\infty}(x)\,\om)  
=-\tfrac{3}{4}\,(I+3 P_{\om})\,\frac{v}{r^2}-\tfrac{3}{4}\,(I-3 P_{\om})\,\frac{v}{r^4}\,.
\end{equation}

Using the obvious scaling relation
$$
\phi_{\var}[v](x)=\Phi_{1,r_{\var}/\var}[v](x/\var)\,,
$$
we see that (\ref{pipo}) and (\ref{pippo}) become, for $r \in [\var, r_{\var}]$
\begin{equation}
\label{phivar}
\begin{aligned}
{\phi}_{\var}[v](x)&=-\,\bigg[4\,\a_1(\var)\,r^2+2\,\b_1(\var)
+\frac{\g_1(\var)}{r}-\frac{\de_1(\var)}{r^3}\bigg]\,(I-P_{\om})\,v 
\\
&-\,2\,\bigg[\a_1(\var)\,r^2+\b_1(\var)+\frac{\g_1(\var)}{r} 
+\frac{\de_1(\var)}{r^3}\bigg]\,P_{\om}\,v
\\
&= A(r)\,(I-P_{\om})\,v\,+B(r)\,P_{\om}\,v
\end{aligned}
\end{equation}
with 
\begin{equation}
\label{as1}
\a_1(\var)=-\tfrac38+O(\var^{2/3})\,,
	\qquad\b_1(\var)=\tfrac98\,\var^{2/3}+O(\var^{4/3}), 
\end{equation}
while
\begin{equation}
\label{as2}
\g_1(\var)=-\tfrac34\,\var+O(\var^{5/3})\,, 
\qquad\de_1(\var)=\tfrac14\,\var^{3}+O(\var^{11/3})\,. 
\end{equation}

Furthermore, for $r \in [\var, r_{\var}]$
\begin{equation} 
\label{gradvar}
\begin{aligned} 
\grad{\phi}_{\var}[v]&=-\,(a(r)+b(r))\,\om\otimes(I-P_{\om})\,v
\\
&+\,b(r)\,\bigg[(I-P_{\om})\,v\otimes\om
	+v\cdot\om\,(I-3\om\otimes\om)\,\bigg]
\end{aligned} 
\end{equation}
with
\begin{eqnarray}  
\label{ab}
a(r)&=&6\,\left(\a_1\,r+\frac{\de_1}{r^4}\right)\,,
\\
\label{ab2}
b(r)&=&2\,\a_1\,r-\frac{\g_1}{r^2}-3\,\frac{\de_1}{r^4}\,. 
\end{eqnarray}
We also record the following formulas for scalar products~:
\begin{equation} 
\label{coso1}
\begin{aligned} 
\grad\phi_{\var}[v_k]:\grad\phi_{\var}[w(x_k)]=
[(a+b)^2+b^2]\,(v_k\cdot w(x_k)-P_{\om}v_k\cdot P_{\om}w(x_k))
\\
+\,6\,b^2\,P_{\om}v_k \cdot P_{\om}w(x_k)\,,
\end{aligned}
\end{equation}
and
\begin{equation}
\label{coso2}
\begin{aligned}
\grad\phi_{\var}[v_k]:\grad w(\cdot+x_k)&=
-\,(a+b)\,\om\cdot\grad(v_k\cdot w)+\,b\,\om\cdot(v_k\cdot\grad w)
\\
&+\,v_k\cdot\om\,(a-3b)\om\cdot(\om\cdot\grad w)
	+\,b\,v_k\cdot\om\,\nabla\cdot w\,.
\end{aligned}
\end{equation}
In the last formula, we have kept the term $\nabla\cdot w$, although all the 
vector fields $w$ considered in this work are solenoidal.

\section{Passing to the limit} \label{esti}

First, we prove that the correctors defined in (\ref{Axp}) converge weakly to
$0$ in $H^1_0(\Om)$ in the vanishing $\var$ limit.

\subsection{Weak convergence of $\Cal{A}_{\var}$} \label{cdc}

Observe that
\begin{equation}
\label{integrale}
\int_{\var\leq|z|\leq r_{\var}}G(|z|)\,P_{z/|z|}v_k\cdot P_{z/|z|}w(x_k)dz
=\tfrac{4\pi}{3}\,v_k\cdot w(x_k)\,\int_{\var\leq r\leq r_{\var}}G(r)\,r^2\,dr
\end{equation}
for each function $G$ for which the integral on the right-hand side makes 
sense. Therefore, using (\ref{phivar}) and (\ref{integrale}), we obtain
\begin{equation}
\|\phi_{\var}[v]\|^2_{L^2(\Om)}=
\tfrac{4\pi}{3}\,|v|^2\,\left(\int^{r_{\var}}_{\var}r^2 \,(2A^2+B^2)\,dr+\var^3\right)
\end{equation}
where the last term comes from the integral on $B_{0,\var}$.

Since 
$$
\begin{aligned} 
\int^{r_{\var}}_{\var} r^2\,(2 A^2 + B^2)\,dr
&\leq
C\,\bigg((\a_1)^2\,(r_{\var}^7-\var^7)+|\a_1\,\b_1|(r_{\var}^5 -\var^{5})
\\
&+\,|\a_1\,\g_1|\,(r_{\var}^4-\var^4)+(\b_1)^2\,(r_{\var}^3-\var^3)
\\
&+\big(|\b_1\,\g_1|+|\a_1\,\de_1|\big)\,(r_{\var}^2-\var^2)
\\
&+\,|\b_1\,\g_1|\,\log (r_{\var}/\var)+(\g_1)^2\,(r_{\var}-\var)
\\
&+|\de_1\,\g_1|\,\left(\frac{1}{\var}-\frac{1}{r_{\var}}\right)
+(\de_1)^2\,\left(\frac{1}{\var^3}-\frac{1}{r_{\var}^3}\right)\bigg), 
\end{aligned}
$$
we obtain from (\ref{as1})-(\ref{as2})~:
\begin{equation}
\label{normal2}
\begin{aligned}
{}&\left\|\sum_{k=1}^{N}\phi_{\var}[v_k]\right\|^2_{L^2 (\Om)}
=\sum_{k=1}^{N}\|\phi_{\var}[v_k]\|^2_{L^2 (\Om)}=
\\
&\tfrac{4\pi}{3}\,\frac{1}{N}\sum_{k=1}^{N}|v_k|^2\,
\bigg(O(N\var^3)+O(N\var^{7/3}) 
+O(N\var^{9/3})+O(N\var^{11/3}\,|\log\var|)\bigg)
\\
&\qquad\qquad\qquad\qquad\qquad\qquad\qquad\qquad\qquad
\qquad\qquad\qquad\qquad\qquad\to 0\,,
\end{aligned}
\end{equation}
so that $\Cal{A}_{\var}\to 0$ in $L^2(\Om)$.

Next, we consider (\ref{coso1}) with $w(x_k)=v_k$. Since
\begin{equation}
\label{coso3}
\grad\phi_{\var}[v_k]:\grad\phi_{\var}[v_k]
=[(a+b)^2+b^2]\,|v_k|^2+ [(a+b)^2+ 5\,b^2] \,(P_{\om}v_k)^2\,,  
\end{equation}
we obtain
$$
\begin{aligned}
{}&\|\grad\phi_{\var}[v_k]\|^2_{L^{2}(B_{0,r_{\var}}\setminus B_{0,\var})}
=
\tfrac{16\pi}{3}\,|v_k|^2\int_{\var}^{r_{\var}}[(a+b)^2+ 2\,b^2]\,r^2\,dr
\\
&\quad=\tfrac{16\pi}{3}\,|v_k|^2\,\bigg(\tfrac{72}{5}\,\a_1^2\,
(r_{\var}^5-\var^5)-3\,\g_1^2\,(r_{\var}^{-1}-\var^{-1})
	-\tfrac{27}{5}\,\de_{1}^2 \,(r_{\var}^{-5}-\var^{-5})
\\
&\quad\quad-\,12\,\a_1\g_1\,(r_{\var}^2-\var^2)
	+24\,\a_1\de_1\,\ln\frac{r_{\var}}{\var}
		-2\,\g_1\de_1\,(r_{\var}^{-3}-\var^{-3})\bigg)
\end{aligned}
$$
which, together with (\ref{as1})-(\ref{as2}), gives
\begin{equation}
\label{eqe}
\|\grad \phi_{\var}[v_k]\|^2_{L^{2}(B_{0,r_\var}\setminus B_{0,\var})}
	\leq C\,\var\,|v_k|^2\,. 
\end{equation}
Hence, setting $r_\var=\var^{1/3}$,
$$
\begin{aligned}
\|\grad \Cal{A}_{\var}\|^2_{L^{2}(\Om)}
&\leq C\,\var\sum_{k=1}^N \,|v_k|^2
\\
&=C\,\bigg(\int_{\Om}\int_{v\in\R^3}F_{N}(x,v)\,|v|^2\,dv\bigg)\leq C'\,. 
\end{aligned}
$$
Therefore, there exists a subsequence such that $\Cal{A}_{\var}\deb\Cal{A}$ 
in $H^1_0(\Om)$. Since $\Cal{A}_{\var}\to 0$ in $L^2(\Om)$, $\Cal{A}=0$ and 
the whole sequence $\Cal{A}_{\var}\deb 0$ in $H^1_0(\Om)$.

\subsection{Weak convergence of $\Cal{B}_{\var}$} \label{cdcb}

Next we prove that the sequence of correctors given by definition (\ref{Bexp})  
converges weakly to $0$ in $H^1_0(\Om)$.

First we estimate $\|\Cal{B}_{\var}[w]\|_{L^2(\Om)}$ and
$\|\grad\Cal{B}_{\var}[w]\|_{L^2(\Om)}$. In order to do so, we consider the 
solution $\Psi_{s,r}[\phi]$ of Stokes problem (\ref{probcor3}) with $w = \phi$.

Using $\Psi_{\var, 2\var} - \Psi_{\var, \var^{1/3}}$ as test function in problem (\ref{probcor3}) with $s=\var$, $r= \var^{1/3}$, we see that
\begin{equation}\label{pc}
\|\grad \Psi_{\var, \var^{1/3}}\|_{L^{2}(B_{0,r_\var}\setminus B_{0,\var})}
\leq\|\grad \Psi_{\var, 2\var}\|_{L^{2}(B_{0, 2\,\var} \setminus B_{0,\var})}.
\end{equation}
According to \cite{allaire}, Lemma 2.2.5, formula 2.2.37, p.240, the following
holds for each $\eta \in ]0,1[$ and $u\in H^1(B_{0,1})$)
\begin{equation}\label{ae}
\|\grad \Psi_{\eta, 1}[u]\|_{L^2(B_{0,1}\setminus B_{0,\eta})}\leq C \, 
\bigg(\|\grad u\|_{L^2(B_{0,1})}+ \eta^{1/2}\|u\|_{L^2(B_{0,1})} \bigg),
\end{equation}
where the constant $C$ is uniform in $\eta$ and $u$.
Observe that $u(x) = \phi(2\var\,x)$ satisfies
\begin{equation}\label{cci}
\|u\|_{L^2(B_{0,1})}= (2\var)^{-3/2}\,\|\phi\|_{L^2(B_{0, 2\var})},
\end{equation}
\begin{equation}\label{bbi}
\|\grad u\|_{L^2(B_{0,1})}= (2\var)^{-1/2}\,\|\grad \phi\|_{L^2(B_{0, 2\var} )},
\end{equation}
and
\begin{equation} \label{ii}
\|\grad \Psi_{\var, 2\var}[\phi]\|_{L^2(B_{0,1\var}\setminus B_{0,\var})}
= (2\var)^{1/2}\|\grad \Psi_{\eta, 1}[u]\|_{L^2(B_{0,1}\setminus B_{0,1/2})}.
\end{equation}
Using successively (\ref{ii}), (\ref{pc}), (\ref{ae}) and (\ref{cci}), (\ref{bbi}), we 
see that
$$ 
\|\grad \Psi_{\var,  \var^{1/3}}[\phi]\|_{L^{2}(B_{0, \var^{1/3}}\setminus B_{0,\var})}
\leq (2\var)^{1/2}\,\|\grad \Psi_{\eta, 1}[u]\|_{L^2(B_{0,1} \setminus B_{0, 1/2})} 
$$
$$ 
\leq C\,\bigg((2\var)^{1/2}\,\|\grad u\|_{L^2(B_{0,1})}
+ 
(1/2)^{1/2}\,(2\var)^{1/2}\,\|u\|_{L^2(B_{0,1})}\bigg) 
$$
\begin{equation} \label{affr}
= C \,\bigg(\|\grad \phi\|_{L^2(B_{0, 2\var} )}
+ (1/2)^{1/2}\,(2\var)^{-1}\,\|\phi\|_{L^2(B_{0, 2\var} )}\bigg)\,. 
\end{equation}
Assuming that $\phi$ is smooth and $\phi(0)=0$ implies that
\begin{equation}\label{stimstab}
\|\grad \psi_{\var}[\phi]\|_{L^{2}(B_{0, \var^{1/3}}\setminus B_{0,\var})} 
\le \hbox{Const.}\var^{3/2}. 
 \end{equation}
Since
$$
\|\grad \Cal{B}_{\var}[w]\|_{L^2 (\Om)}\leq 
\sum_{k=1}^N\|\grad \phi_{\var}[w(x_k)]\|_{L^2}
+\sum_{k=1}^N || \nabla \psi_{\var}[w(\cdot+x_k)-w(x_k)]||_{L^2}\,,
$$
it follows from (\ref{eqe}) with $r_{\var}=\var^{\frac 13}$ and $v_k=w(x_k)$ 
and (\ref{stimstab}) that 
$$
\|\grad\Cal{B}_{\var}[w]\|_{L^2(\Om)} <\hbox{Const.}
$$
Hence there is a subsequence s.t. $\mathcal{B}_{\var}\deb\mathcal{B}$ in $H^1_0(\Om)$.

On the other hand
$$
\begin{aligned}
\sum_{k=1}^{N}\psi_{\var}[w(\cdot+x_k)](x-x_k)
	=\sum_{k=1}^{N}\phi_{\var}[w(x_k)](x-x_k) 
\\
+\sum_{k=1}^{N}\psi_{\var}[w(\cdot+x_k)-w(x_k)](x-x_k)\,.
\end{aligned}
$$ 
By Poincar\'e's inequality and (\ref{stimstab}), 
$$ 
\begin{aligned}
\sum_{k=1}^{N} ||\psi_{\var}[w(\cdot+x_k)-w(x_k)](x-x_k)||_{L^2}^2 
\\
\le \var^{2/3}\,  
\sum_{k=1}^{N} ||\nabla \psi_{\var}[w(\cdot+x_k)-w(x_k)](x-x_k)||_{L^2}^2 
\\
\le Cst\, N\, \var^{11/3} .
\end{aligned}
$$
Using (\ref{normal2}) with $v_k=w(x_k)$ shows that 
$$
\mathcal{B}_{\var}=\sum_{k=1}^{N}\psi_{\var}[w(\cdot+x_k)]\to 0
$$ 
in $(L^2(\Om))^3$, so that $\mathcal{B}=0$ and the whole sequence 
$$
\mathcal{B}_{\var}\deb 0\hbox{ in }(H^1_0(\Om))^3\,.
$$
Finally, notice that 
$$
\|\sum_{k=1}^{N}\psi_{\var}[w(\cdot+x_k)-w(x_k)]\|_{(H^1_0(\Om))^3}\to 0\,.
$$

\subsection{Limit of $\int_{\Om} \grad \Cal{A}_{\var}: \grad W_{\var}$} \label{kk1}

We have
\begin{equation}
\label{stab2}
\begin{aligned}
\int_{\Om} \grad \mathcal{A}_{\var}: \grad W_{\var}dx
=\int_{\Om}\grad \mathcal{A}_{\var}: \grad w dx
-\int_{\Om}\grad \mathcal{A}_{\var}: \grad \Cal{B}_{\var}[w]dx
\\ 
=\int_{\Om}\grad \mathcal{A}_{\var}: \grad w dx
-\int_{\Om}\grad \mathcal{A}_{\var}: 
\grad (\sum_{k=1}^{N}\psi_{\var}[w(\cdot+x_k)-w(x_k)])dx
\\
-\sum_{k=1}^N\int_{\var\leq|z|\leq r_{\var}}
	\grad\phi_{\var}[v_k]:\grad\psi_{\var}[w(x_k)]dz 
\end{aligned}
\end{equation}
In view of (\ref{A}), we see that
$$
\lim_{\var\to 0}\int_{\Om} \grad \mathcal{A}_{\var}: \grad wdx=0\,.
$$
Moreover, since 
$$
\grad (\sum_{k=1}^{N}\psi_{\var}[w(\cdot+x_k)-w(x_k)])\to 0
\hbox{ in }(L^2(\Om))^9
$$
we also have
$$
\int_{\Om}\grad \Cal{A}_{\var}: 
	\grad (\sum_{k=1}^{N}\psi_{\var}[w(\cdot+x_k)-w(x_k)])dx\to 0\,.
$$

Next we estimate (\ref{stab2}). Recalling (\ref{integrale}), we have
\begin{equation}
\int_{\var\leq|z|\leq r_{\var}}\grad\phi_{\var}[v_k]:\grad\phi_{\var}[w(x_k)]dx
= 
4\pi \, v_k \cdot w(x_k) \, 
	\bigg(\int_{\var\leq r\leq r_{\var}} \mathcal{F}(r)\,r^2\,dr \,\bigg),
\end{equation}
where (according to (\ref{coso1}))
\begin{equation}
\mathcal{F}(r)= \frac23\,[(a+b)^2+b^2]+ 2\,b^2.
\end{equation}
Therefore,
$$
\begin{aligned}
\int_{\var}^{r_{\var}} r^2\, \mathcal{F}(r)\,dr
=
\tfrac{32}3\, \a_1^2 \,(r_{\var}^5-\var^5)- 6\, \de_1^2\,(r_{\var}^{-5}-\var^{-5}) 
- \tfrac{10}3\,\g_1^2 \,(\frac{1}{r_{\var}}-\frac{1}{\var})
\\
- \tfrac{32}3\, \a_1\,\g_1\,(r_{\var}^2-\var^2)
	-4 \,\g_1\,\de_1\,(r_{\var}^{-3}-\var^{-3})
\\
=\tfrac32\,\var + O(\var^{5/3}) . 
\end{aligned}
$$
Finally,
\begin{equation}
\int_{\var\leq|z|\leq r_{\var}}\grad\phi_{\var}[v_k]:\grad\phi_{\var}[w(x_k)]dz
= 6\pi\, \var v_k\, \cdot w(x_k)+O(\var^{5/3}), 
\end{equation}
and
$$ 
\begin{aligned}
\lim_{\var\rightarrow 0}&\sum_{k=1}^N \int_{\var\leq|z|\leq r_{\var}}
\grad\phi_{\var}[v_k]:\grad\phi_{\var}[w(x_k)]dx
\\
&=\lim_{\var\rightarrow 0}
\sum_{k=1}^{N} \bigg(6\pi\, \var \,v_k \cdot w(x_k)    + O(\var^{5/3})\bigg) 
\\
&=\lim_{\var\rightarrow 0}\,
(6\pi + O(\var^{2/3}))\, (\var N) 
	\int_{\Om}\int_{\R^3}   v\cdot w(x)\, F_N(x,v)\, dvdx 
\\
&\qquad\qquad= 6\pi \, \int_{\Om} j(x)\cdot w(x)\, dx\,.
\end{aligned}
$$
Therefore, we conclude that
\begin{equation}
\lim_{\var\rightarrow 0}\int_{\Om} \grad \Cal{A}_{\var}: \grad W_{\var}dx
	= 6\pi \, \int_{\Om} j(x)\cdot w(x)\, dx\,.
\end{equation}

\subsection{Limit for $\int_{\Om} \grad U_{\var}:\grad \Cal{B}_{\var}$} \label{kk2}

We have
\begin{equation}
\begin{aligned}
\label{sstab1}
\int_{\Om} \grad U_{\var}:\grad \mathcal{B}_{\var}[w]dx
&=
\int_{\Om} \grad U_{\var}:
\grad (\sum_{k=1}^{N}\psi_{\var}[w(\cdot+x_k)-w(x_k)])dx
\\
&+\sum_{k=1}^N\int_{\var\leq|z|\leq r_{\var}}
	\grad U_{\var}(\cdot + x_k):\grad\phi_{\var}[w(x_k)]dz
\end{aligned}
\end{equation}
Since
$$
\grad (\sum_{k=1}^{N}\psi_{\var}[w(\cdot+x_k)-w(x_k)])\to 0
$$ 
in $(L^2(\Om))^9$, we see that
$$
\int_{\Om} \grad U_{\var}:
\grad (\sum_{k=1}^{N}\psi_{\var}[w(\cdot+x_k)-w(x_k)])dx\to 0\,.
$$
In order to estimate (\ref{sstab1}), we first integrate by parts, denoting by
$n$ the outward unit normal vector to the sphere $\pa B_{x_k,r_{\var}}$:
$$
\begin{aligned}
I_k&=\int_{\var\leq|z|\leq r_{\var}}
	\grad U_{\var}(\cdot + x_k):\grad\phi_{\var}[w(x_k)]dz
\\
&=-\int_{\var\leq |z| \leq r_{\var}}
	U_{\var}(\cdot + x_k)\cdot\Lap\phi_{\var}[w(x_k)]dz
\\
&+\int_{\pa B_{x_k,r_{\var}}} 
	(n\cdot \grad\phi_{\var}[w(x_k)](\cdot - x_k))\cdot U_{\var}dz\,.
\end{aligned}
$$
Next we use the definition of $\phi_{\var}$  to compute
\begin{equation}
\label{for}
I_k
=\int_{\pa B_{x_k,r_{\var}}}\bigg(n\cdot\grad\phi_{\var}[w(x_k)](\cdot -x_k)
-\pi_{\var}(\cdot - x_k) \,n\bigg)\cdot U_{\var}dz\,.
\end{equation}
At this point, we observe that  $n=\omega = \frac{x-x_k}{|x-x_k|}$ (that
is, $\omega$ is ``centered'' on $x_k$ instead on the origin as in section 
\ref{soluz}). Since 
$\grad\phi_{\var}(x) = \frac1{\var}\,\grad \Phi_{1,r_{\var}/\var}(x/\var)$ 
it follows from (\ref{oubli}) that
$$
\begin{aligned} 
I_k= \frac1{\var}\,\int_{\pa B_{x_k,r_{\var}}} 
	\bigg(\omega\cdot \grad \Phi_{1,r_{\var}/\var}[w(x_k)](\tfrac{x - x_k}\var) 
- \Pi_{1,\frac{r_{\var}}{\var}}(\tfrac{x - x_k}\var\,\om \bigg) \cdot U_{\var}dx
\\
= \frac1{\var}\,\int_{\pa B_{x_k,r_{\var}}} \bigg(\omega\cdot  
	\grad \Phi_{1,\infty}[w(x_k)](\tfrac{x - x_k}\var) 
		-\Pi_{1,\infty}(\tfrac{x - x_k}\var)\,\om
\\
+\,3\,\left(\frac{\var}{r_{\var}}\right)^2 \,(I-3P_{\om})\,w(x_k)
	+\,O\left(\frac{\var}{r_{\var}}\right)^3 \bigg) \cdot U_{\var}dx\,.
\end{aligned}
$$
By (\ref{forza}), 
$$ 
\begin{aligned}
I_k = \frac{\var}{r_{\var}^2}\int_{\pa B_{x_k,r_{\var}}}\!\!\!
\bigg( - \tfrac34\,(I+3P_{\om})\,w(x_k) +3\,(I-3P_{\om})\,w(x_k)
\\
-\tfrac34\,(I-3P_{\om})\,w(x_k)\,\frac{\var^2}{r_{\var}^2}\bigg)
	\cdot U_{\var}dx+O\left(\frac{\var^2}{r_{\var}^3}\right)\,, 
\end{aligned}
$$
so that
$$ 
I_k\!\!=\!\!\frac{\var}{r_{\var}^2}\int_{\pa B_{x_k,r_{\var}}}\!\!\!
\bigg(- \tfrac34\,(I+3P_{\om})\,w(x_k) +\,3\,(I-3P_{\om})\,w(x_k)\bigg)
\cdot U_{\var}dx+O\left(\frac{\var^2}{r_{\var}^3}\right).
$$
Notice that the same result is obtained in (\cite{allaire})) by a somewhat
different procedure.

At this point, we claim the following strong limits in $(H^{-1}(\R^3))^3$ that 
hold for any $G\in (C_b(\R^3))^3$ --- for a proof, see sec. (\ref{la})) in the 
appendix below:
\begin{equation}
\begin{aligned}
\label{uffa}
\sum_{i=1}^{N}r_{\var}\,G(x_k)\,\de_{\pa B_{x_k, r_{\var}}}
	\to&4\pi\,\rho(x)\,G(x)&&\hbox{ in }(H^{-1}(\R^3))^3\,,
\\
\sum_{i=1}^{N}r_{\var}\,G(x_k)\cdot\om \,\om \,\de_{\pa B_{x_k,r_{\var}}}
\to &\tfrac{4\pi}3\,\rho(x)\, G(x) &&\hbox{ in }(H^{-1}(\R^3))^3\,.
\end{aligned}
\end{equation}
Since $U_{\var}\to U$ converges weakly in $(H^1(\Omega))^3$, we get
\begin{equation}
\sum_{k=1}^N 
\int_{\var\leq|z|\leq r_{\var}}\grad U_{\var}:\grad\psi_{\var}[w(\cdot+x_k)]dz
\to-\,6\pi\,\int_{\Om}w(x)\cdot U \,\rho(x)dx\,,
\end{equation}
so that
\begin{equation}
\lim_{\var\to 0}
\int_{\Om}\grad U_{\var}:\grad \mathcal{B}_{\var}dx
=-\,6\pi\,\int_{\Om} w(x)\cdot U \,\rho(x)dx.
\end{equation}

\subsection{The limit  equation}

We start from the weak formulations established in section \ref{wee} (that is, 
equations (\ref{fdstokes}) and (\ref{fdnavstok})).

\subsubsection{The Stokes case}

In view of the results established in sections \ref{kk1} and \ref{kk2}, we pass 
to the limit in the Stokes problem (\ref{prob})-(\ref{bc})
$$
\int_{\Omega}\grad U\cdot\grad{w}dx=
\int_{\Omega}{g}\cdot wdx +
6\pi \left(\int v\cdot w fdxdv-\int w\cdot U \rho dx\right)\,.
$$
Therefore, $U$ is the a weak solution of 
\begin{equation}
\left\{
\begin{array}{ll}
-\Lap U -\,6\pi \,(j- \rho\, U)+ \grad \Pi= g\,,
\\
\nabla\cdot U=0\,,
\\
U |_{\pa \Om}=0\,.
\end{array}
\right.
\end{equation}
Since the problem above has at most one weak solution, the whole sequence
$U_{\var}$ converges to $U$ in $(H^1(\Omega))^3$.

In addition
$$
\bar{u}_{\var}=U_{\var}+\sum_{k=1}^{N}\phi_{\var}[v_k](x-x_k)\to U
$$
in $(L^2(\Om))^3$, as can be seen from (\ref{normal2}).

This finishes the proof of theorem 1, assuming (\ref{uffa}) --- whose proof 
is deferred to the appendix below

\subsubsection{The Navier-Stokes case}

Likewise, for the Navier--Stokes problem (\ref{probNS}), (\ref{bc}) in the 
limit as $\var\to 0$
$$
\begin{aligned}
\nu\int_{\Omega}\grad U \cdot\grad{w}dx
&=\int_{\Omega} U\otimes U :\grad wdx
\\
&+\int_{\Omega}{g}\cdot wdx
+6\pi 
\left(\int v\cdot wfdxdv-\int w\cdot U\rho dx\right)\,.
\end{aligned}
$$

Given $\rho$, $j$ and $g$, there exists $\nu_0>0$ large enough, so that, 
for each $\nu > \nu_{0}$, the problem
\begin{equation}
\left\{
\begin{array}{ll}
U\cdot\grad U-\nu\Lap U(x) -\,6\pi \,(j- \rho\, U)+ \grad \Pi= g\,,
\\
\nabla\cdot U=0\,,
\\
U |_{\pa \Om}=0\,,
\end{array}
\right.
\end{equation}
has a unique weak solution $U\in (H^1_0(\Om))^3$.

Hence the whole sequence $U_{\var}$ converges weakly to $U$ in 
$(H^1_0(\Om))^3$ as $\var\to 0$.

As in the Stokes case, (\ref{normal2}) implies that 
$$
\bar{u}_{\var}=U_{\var}+\sum_{k=1}^{N}\phi_{\var}[v_k](x-x_k)\to U
$$
in $(L^2(\Om))^3$. This completes the proof of theorem 2 --- assuming 
again that the limits in(\ref{uffa}) hold.

\section{Appendix}\label{append1}

\subsection{Proof of (\ref{uffa})}\label{la}

We closely follow the method described in \cite{cioramurat} in the periodic
setting.
Given  $G\in (C_{b}(\R^3))^3$, we consider two  auxiliary problems 
in $\bigcup_{k=1}^{N}B_{x_k, r_{\var}}$~: for $k=1,\ldots,N$
\begin{eqnarray*}
\left\{ \begin{array}{ll}
-\Lap \xi_{\var}=-3 \,G(x_k),& 
\\
\frac{\pa \xi_{\var}}{\pa n}\mid_{\pa B_{x_k,r_{\var}}}
	= r_{\var} \,G(x_k),&
\\
\end{array}
\right.
\end{eqnarray*}
and
\begin{eqnarray*}
\left\{ \begin{array}{ll}
-\!\!\Lap\! \chi_{\var}\!=-G(x_k)-\frac{r}{r_{\var}}\,\bigg(6\,(G(x_k)\cdot\om)\, \om 
+2\,G(x_k)\bigg)\!\!+\,3 \,G(x_k), &
\\
\frac{\pa \chi_{\var}}{\pa n}\mid_{\pa B_{x_k,r_{\var}}}
	= r_{\var}\,  (G(x_k)\cdot \om )\,\om\,. &
\\
\end{array}
\right.
\end{eqnarray*}
Next we extend $\xi_{\var}$ and $\chi_{\var}$ by $0$ in the complement
of  $\bigcup_{k=1}^{N}B_{x_k, r_{\var}}$. 
Computing the Laplacian of $\xi_{\var}$ and $\chi_{\var}$ in the sense of 
distributions in the whole Euclidean space, we get
\begin{eqnarray}
-\Lap \xi_{\var}&=&\!\!\! -3 \sum_{k=1}^{N} G(x_k)\, 1_{ B_{x_k,r_{\var}}}
+\sum_{k=1}^{N} r_{\var} \,G(x_k) \,\de_{\pa B_{x_k,r_{\var}}}\nonumber
\\
&=&\!\!\!-3 \,N \,1_{ B_{0,r_{\var}}}\ast (G\rho_N)
+\sum_{k=1}^{N} r_{\var} \, G(x_k) \,\de_{\pa B_{x_k,r_{\var}}}, \label{secmu1}
\end{eqnarray}
\begin{eqnarray}
-\Lap \chi_{\var}
&=&\!\!\!  -\sum_{k=1}^{N} G(x_k) \,1_{ B_{x_k,r_{\var}}}\nonumber
\\
&+&
\sum_{k=1}^{N}\bigg\{\frac{r}{r_{\var}}\,\bigg(6\, (G(x_k)\cdot\om)\, \om 
+ 2 \,G(x_k)\bigg)-3 \,G(x_k)\bigg\} \,1_{ B_{x_k,r_{\var}}}\nonumber
\\
&+&\sum_{k=1}^{N} r_{\var} \,(G(x_k)\cdot\om)\, \om 	
	\,\de_{\pa B_{x_k,r_{\var}}}\label{secmu2}
\\
&=&\!\!\!- N \,1_{ B_{0,r_{\var}}}\ast (G\rho_N)
+\sum_{k=1}^{N} r_{\var} \, (G(x_k)\cdot\om)\, \om 
\,\de_{\pa B_{x_k,r_{\var}}}\nonumber
\\
&+&
\sum_{k=1}^{N} \bigg\{\frac{r}{r_{\var}}\,\bigg(6\,(G(x_k)\cdot\om)\, \om 
+
2 \,G(x_k)\bigg)-3 \,G(x_k)\bigg\} \,1_{ B_{x_k, r_{\var}}} \,.\nonumber
\end{eqnarray}
The solutions of the two auxiliary problems above are
\begin{eqnarray*}
\xi_{\var}(x)&=&\sum_{i=1}^{N} 
\frac{|x-x_k|^2-r_{\var}^2}{2} \,G(x_k)\,1_{ B_{x_k,r_{\var}}},
\\
\chi_{\var}(x)&=&\sum_{i=1}^{N} 
\bigg(\frac{|x-x_k|^3}{r_{\var}}-|x-x_k|^2\bigg)\, (G(x_k)\cdot \om )\,\om \, 
1_{ B_{x_k,r_{\var}}},
\end{eqnarray*}
while their gradients are given by
\begin{eqnarray*}
\grad\xi_{\var}(x)&=&
\sum_{k=1}^{N} |x-x_k|\,(\om\otimes G(x_k))\,1_{ B_{x_k,r_{\var}}},
\\
\grad\chi_{\var}(x)&=&
\sum_{k=1}^{N} \frac{|x-x_k|^2}{r_{\var}} \,
((\om\otimes G(x_k))\cdot\om) \,\om \,1_{ B_{x_k,r_{\var}}}
\\
&+&\sum_{k=1}^{N}\bigg(\frac{|x-x_k|^2}{r_{\var}}-|x-x_k|\bigg)\,
\\
&\times&\bigg((G(x_k)\cdot \om)\,I+G(x_k)\otimes\om\bigg)\,1_{ B_{x_k,r_{\var}}}.
\end{eqnarray*}
with $\omega = \frac{x-x_k}{|x-x_k|}$.

Then, we estimate
\begin{eqnarray*}
\|\xi_{\var}\|_{L^2(\R^3)}^2
	&\leq&\hbox{Const.}\,r_{\var}^4\,\int G^2 \rho_Ndx=O(r_{\var}^4),
\\
\|\grad\xi_{\var}\|_{L^2(\R^3)}^2
	&\leq&\hbox{Const.}\, r_{\var}^2 \,\int G^2 \rho_Ndx=O(r_{\var}^2),
\\
\|\chi_{\var}\|_{L^2(\R^3)}^2
	&\leq&\hbox{Const.}, r_{\var}^4\,\int G^2 \rho_Ndx=O(r_{\var}^4),
\\
\|\grad\chi_{\var}\|_{L^2(\R^3)}^2
	&\leq&\hbox{Const.}t\, r_{\var}^2 \,\int G^2 \rho_Ndx=O(r_{\var}^2).
\end{eqnarray*}
Therefore, $\xi_{\var}$ and $\chi_{\var}\to 0$ in $(H^{1})^3$, so that both
\begin{equation}
\label{lim0}
\Delta \xi_{\var}\hbox{ and }\Delta \chi_{\var}\to 0\hbox{ in }(H^{-1})^3\,.
\end{equation}

Next, we recall that $N \,1_{ B_{0,r_{\var}}}\deb\tfrac{4\pi}3\de_0$ weakly 
in the sense of measures; hence
\begin{equation}
\label{lim1}
N \,1_{ B_{0, r_{\var}}}\ast (G\rho_N)\deb\tfrac{4\pi}3\,\rho \,G
\end{equation}
weakly in the sense of measures. Furthermore
$$
\|N \,1_{ B_{x_k, r_{\var}}}\ast (G\rho_N)\|_{L^\infty}\leq \|G\|_{L^\infty}
$$
so that, by the Rellich compactness theorem, the limit (\ref{lim1}) holds in
the strong topology of $(H^{-1}_{loc})^3$. 

Going back to (\ref{secmu1}), we conclude from (\ref{lim0}) that
$$
\sum_{k=1}^{N}r_{\var}\,G(x_k)\,\de_{\pa B_{x_k,r_{\var}}}
	\to 4\pi\,\rho\,G\hbox{ in }(H^{-1}_{loc})^3
$$
strongly.

\smallskip
Next, we apply the same procedure to the second term on the right hand 
side of (\ref{secmu2}).

First, we observe that the last term on that right hand side is bounded 
in $L^\infty(\R^3)$ by $11\|G\|_{L^\infty}$, while
\begin{equation}
\label{lim2} 
\int_{\Om} \bigg(\phi(x)\cdot\sum_{k=1}^{N} \bigg\{\frac{r}{r_{\var}}\,\bigg(6\,(G(x_k)\cdot\om)\om 
	+2G(x_k)\bigg)-3G(x_k)\bigg\}1_{ B_{x_k, r_{\var}}} \bigg)dx \to 0
\end{equation}
for each $\phi\in(\mathcal{D}(\R^3))^3$. Applying the Rellich compactness
theorem again shows that the convergence (\ref{lim2}) holds in the strong
topology of $(H^{-1}_{loc})^3$.

Going back to (\ref{secmu2}) and using (\ref{lim0}), (\ref{lim1}) and (\ref{lim2})
shows that
$$
\sum_{k=1}^{N}r_{\var}(G(x_k)\cdot\om)\om\de_{\pa B_{x_k, r_{\var}}}
\to\tfrac{4\pi}{3 }\,\rho G\hbox{ in }(H^{-1})^3
$$
strongly.

\subsection{Solution of Stokes' problem in an annulus}\label{append}

We first prove the explicit formula for $\Phi_{1,R}$ in (\ref{pipo}), 
(\ref{pippo}) and (\ref{pipppo}), by the same method as in 
\cite{LL} \S 20. By symmetry, we seek $\Phi_{1,R}$ in the form 
$\Phi_{1,R}[v] = \mathrm{curl}\, \mathrm{curl} (f(r)\, v)$ (where
$r=|x|$). Then
$$ 
\bigg(\frac{\pa^2}{\pa r^2} + \frac2r \, \frac{\pa}{\pa r}\bigg)\,  
\bigg(\frac{\pa^2}{\pa r^2} + \frac2r \, \frac{\pa}{\pa r}\bigg) f(r) 
= \hbox{Const.} 
$$
so that
$$ 
f'(r) = \a\, r^3 + \beta\,r + \gamma + \frac{\de}{r^2}. 
$$
Denoting by $P_{\om}$ the orthogonal projection on $\om=x /r$ we 
arrive at formula (\ref{pipo})~:
\begin{eqnarray*}
\Phi_{1,R}[v](x) &=&-\,(f^{''}+\frac{f^{'}}{r})\, (I-P_{\om})\,v-2\,\frac{f^{'}}{r} \,P_{\om}v
\\
&=&-\, \bigg[4\,\a (R) \,r^2 + 2 \,\b (R) -\frac{\de (R)}{r^3}+\frac{\g (R)}{r}\bigg]\,
(I-P_{\om})\,v
\\
&&-\,2\, \bigg[\a (R) \,r^2 + \b (R) +\frac{\de (R)}{r^3}+\frac{\g (R)}{r}\bigg]\,
P_{\om}\,v.
\\
\end{eqnarray*}
Because of the boundary conditions, the constants $\a$,$\b$,$\g$,$\de$ 
in the formula above satisfy the following system of equations~:
\begin{eqnarray*}
3\,(R^5-1)\,\a &+&(R^3-1)\,\b =\frac{1}{2},
\\
5\,(R^3-1)\,\a &+&3\,(R-1)\,\b =\frac{3}{2},
\\
\tfrac{3}{2}\,\a &+&\tfrac{1}{2}\,\b +\tfrac{1}{4}=\de,
\\
-\tfrac{5}{2}\,\a &-&\tfrac{3}{2}\,\b -\tfrac{3}{4}=\g, 
\end{eqnarray*}
leading to the estimates (\ref{pippo}) and (\ref{pipppo})~:
$$
\a=-\frac{3}{8R^3}+O_{R \to +\infty} (1/R^4), 
\qquad \b=\frac{9}{8R}+O_{R \to +\infty}(1/R^2), 
$$
$$
\g=-\frac34 +O_{R \to +\infty}(1/R), 
\qquad \de=\frac14 +O_{R \to +\infty}(1/R). 
$$

Next we compute the pressure for the above flow; for simplicity, we first 
write down the following  table~:
\begin{eqnarray*}
\label{csempl}
\begin{array}{ll}
\grad r=\om, & \grad \om = \frac{I-\om\otimes\om}{r},\\
\grad\cdot \om=\frac{2}{r},&\om\cdot\grad \om=0,\\
\grad (\om\cdot v) = \frac{v-P_{\om}v}{r},\\
\om\cdot P_{\om}v = \om\cdot v,&\grad\cdot (a\otimes b)
=(\grad\cdot a) \,b + a\cdot\grad b,
\end{array}
\end{eqnarray*}
so that
\begin{eqnarray*}
\grad P_{\om}\,v&=&\frac{v-P_{\om}v}{r}\otimes\om+(\om\cdot v) \,\grad\om
\,\,\,\,\,\mathrm{and}\,\,\,\,\,\om\cdot\grad P_{\om}\,v=0,\\
\grad \cdot P_{\om}v&=&
\frac{2}{r}\,\om\cdot v,\,\,\,\,\,\,\,\,\,\,
\Lap \om=
-\,\frac{2}{r^2}\,\om\, ,\\\
\Lap P_{\om}\,v &=&
2 \,\bigg(\frac{v-3P_{\om}v}{r^2} \bigg) .
\end{eqnarray*}

To find the pressure, observe that
$$
\Lap \Phi_{1,R}[v](x)
=
-\, \bigg[20\,\a -\frac{5\,\a}{r^3}-\frac{3\,\b}{r^3}-\frac{3}{2 \,r^3} \bigg]\,v-\,3\, \bigg[\frac{5\,\a}{r^3}+
\frac{3\,\b}{ r^3}+\frac{3}{2 \,r^3}\bigg] \,P_{\om}\,v
$$
so that, up to some unessential additive constant, 
$$ 
\Pi_{1,R}=- \,\bigg[\bigg(20\,\a \,r- \frac{3}{2 \,r^2}
	-\frac{5\,\a+3\,\b}{r^2}\bigg)\,\om\cdot P_{\om}\,v \bigg].
$$

Letting $R\to +\infty$  in the previous expressions leads to formula (\ref{infinito}) 
(see also \cite{LL})~:
\begin{equation}\label{lab1}
\Phi_{1,\infty}[v](x) =
\tfrac{1}{4}\, (\frac{3}{r}+\frac{1}{r^3})\,(I-P_{\om})\,v
+\tfrac{1}{2}\,(\frac{3}{r}-\frac{1}{r^3})\,P_{\om} \,v , 
\end{equation}
and 
\begin{equation}\label{lab2}
\Pi_{1, \infty}=\frac{3}{2 r^2} \,\om\cdot P_{\om}\,v . 
\end{equation}

Obviously
\begin{eqnarray*}
\Phi_{1,R}[v](x) &=&\Phi_{1,\infty}[v](x) 
	-[4\,\a (R) \,r^2 + 2 \,\b (R)]\,(I-P_{\om})\,v
\\
&-&\,2\,[\a (R) \,r^2 +  \b (R)]\,P_{\om}v+\frac1r\,O\left(\frac{1}{R}\right),
\\
\om\cdot\grad \Phi_{1,R}[v](x)&=&\om\cdot \grad \Phi_{1,\infty}[v](x) 
\\
&-&\,4\,\a(R)\,r \,(2I-P_{\om})\,v+\frac1{r^2}\,O\left(\frac{1}{R}\right),
\end{eqnarray*}
while 
$$ 
\Pi_{1,R}=\Pi_{1,\infty}-20 \,\a(R)\,r \,\om\cdot v+ \frac1{r^2}\,
O\left(\frac{1}{R}\right). 
$$
as $R\to+\infty$.

Formula (\ref{oubli}) follows as a consequence:
$$ 
\begin{aligned}
\om\cdot\nabla \Phi_{1,R}-\Pi_{1,R}\,\om&=
\om\cdot\nabla \Phi_{1,\infty} -\Pi_{1,\infty}\,\om 
\\
&- 8\,\a(R)\,r \,(I-3P_{\om})\,v+\frac1{r^2}\,O\left(\frac{1}{R}\right).  
\end{aligned}
$$
From (\ref{lab1}) and (\ref{lab2}), we arrive at formula (\ref{forza})~:
\begin{equation}
\om\cdot\nabla \Phi_{1,\infty}[v](x) -\Pi_{1,\infty}\,\om
= -\tfrac{3}{4}\,(I+3 P_{\om})\,\frac{v}{r^2}
	-\tfrac{3}{4}\,(I-3 P_{\om})\,\frac{v}{r^4}.
\end{equation}

Finally, we derive formulas (\ref{coso1}), (\ref{coso2}). First rewrite 
formula (\ref{gradvar}) for $\grad{\phi}_{\var}[v]$ with $r \in [\var, r_{\var}]$
in the form
\begin{equation}\label{graphi}
\grad{\phi}_{\var}[v]= -a(r) \,N(v) 
	+ b(r)\,[M(v) + v\cdot \om \,(I-3\, \om\otimes\om )] 
\end{equation}
where $M$ and $N$ are two matrix-valued, linear functions of $v$~:
\begin{eqnarray*}
M(v)&=&(I-P_{\om})\,v\otimes\om -\om \otimes (I-P_{\om})\,v,
\\
N(v)&=&\om \otimes (I-P_{\om})\,v,
\end{eqnarray*}
while
$$
a(r) =6\,(\a_1 \,r+\frac{\de_1}{r^4}),
\qquad
b(r) =2\,\a_1\, r -\frac{\g_1}{r^2} -3\,\frac{\de_1}{r^4}. 
$$
For each $v,w\in\R^3$ and each $\omega\in S^2$
\begin{eqnarray*}
\tfrac{1}{2}M(v):M(w)&=& \!\!\!\!N(v):N(w)= \!\!- M(v):N(w)
	=(v\cdot \!w \!-\!\!P_{\om}v\cdot P_{\om}w),
\\
M(v):\Cal{\om\otimes\om}&=&M(v):I= N(v):\Cal{\om\otimes\om}=N(v):I=0,
\\
\Cal{\om\otimes\om}:I&=&\Cal{\om\otimes\om}:\Cal{\om\otimes\om}=1,
\quad I:I=3. 
\end{eqnarray*}
Now (\ref{coso1}) and  (\ref{coso2}) follow from (\ref{graphi}) by elementary
manipulations involving the identities recalled above.

\bigskip


\begin{thebibliography}{0}


\bibitem{allaire} 
G. Allaire.
\newblock {\em Homogenization of the Navier-Stokes equations in 
open sets perforated with tiny holes}.
\newblock Arch. Rational Mech. Anal., {\bf{113}}, p. 209-259, 1991.

\bibitem{kiva}
A.A.~Amsden, P.J.~O'Rourke, T.D.~Butler.
\newblock{\em A computer program for chemically reactive
flows with sprays}.
\newblock Report \#\,LA-11560-MS, Los Alamos National Laboratory, 1989.

\bibitem{Batch}
G.K. Batchelor.
\newblock{\em Sedimentation in a dilute suspension of spheres}.
J. Fluid Mech. \textbf{52} (1972), 245--268.

\bibitem{baranger}
C. Baranger, L. Desvillettes.
\newblock{\em Coupling Euler and Vlasov equations in the context of sprays: 
the local-in-time, classical solutions}. 
J. Hyperbolic Differ. Equ. \textbf{3} (2006), 1--26. 

\bibitem{CaflLuke}
R.E. Caflisch, J. H.C. Luke.
\newblock{\em Variance in the sedimentation speed of a suspension}.
Phys. Fluids \textbf{28} (1985), 759--760.

\bibitem{CaflRubi}
R.E. Caflisch, J. Rubinstein.
``Lectures on the mathematical theory of multi-phase flows"
Courant Institute Lecture Notes, New-York, 1984. 

\bibitem{cioramurat}
D. Cioranescu, F. Murat.
\newblock {\em Une terme \'etrange venu d'ailleurs}. In
\newblock ``Nonlinear Partial Differential Equations and their Applications", 
Coll\`ege de France Seminar, Vol. 2, Research Notes in Mathematics, 
{\bf{60}}, p.98-138, (1982).

\bibitem{Feuille}
F. Feuillebois.
\newblock{\em Sedimentation in a dispersion with vertical inhomogeneities}.
J. Fluid Mech. \textbf{139} (1984), 145--171.

\bibitem{jabin}
P.-E. Jabin.
\newblock{\em Various levels of models for aerosols}.
Math.Models and Methods in Appl. Sci. \textbf{12} (2002), 903--919.

\bibitem{JabOtto}
P.-E. Jabin, F. Otto.
\newblock{\em Identification of the dilute regime in particle sedimentation}.
Commun. Math. Phys. \textbf{250} (2004), 415--432.

\bibitem{Ladyz}
O.A. Ladyzhenskaya:
``The mathematical theory of viscous incompressible flow". 
Mathematics and its Applications, Vol. 2 Gordon and Breach, 
Science Publishers, New York-London-Paris 1969.
 
\bibitem{LL}
L.D. Landau, E.M. Lifshitz.
\newblock{\em Course of theoretical physics. Vol. 6. Fluid mechanics}. 
Pergamon Press, Oxford, 1987.

\bibitem{hrlv}
V. A. L'vov, E. Ya. Khruslov.
\newblock{\cyr  O vozmushtenii vyazko{\u\i}  neszhimaemo{\u\i}  zhidkosti melkimi
chastitsami}(Russian)
[Perturbation of a viscous incompressible fluid by small particles]
 Theoretical and applied questions of differential equations and algebra
 (Russian) \textbf{267} (1978), 173--177.

\bibitem{orourke}
P.J.~O'Rourke.
\newblock {\em Collective drop effects on vaporizing liquid sprays}.
\newblock PhD thesis, Los Alamos National Laboratory, 1981.

\bibitem{RubiKeller}
J. Rubinstein, J. Keller:
\newblock{\em Particle distribution functions in suspensions}.
Phys. Fluids A \textbf{1} (1989), 1632--1641.

\bibitem{Rubi}
J. Rubinstein:
\newblock{\em  On the macroscopic description of slow viscous flow past 
a random array of spheres}. 
J. Statist. Phys. \textbf{44} (1986), 849--863.






 \end{thebibliography}
\end{document}